\documentclass[12pt, a4paper, reqno]{amsart}

\usepackage{amscd, amssymb, amsthm, enumerate}

\newtheorem{theorem}{Theorem}[section]
\newtheorem{corollary}[theorem]{Corollary}
\newtheorem{lemma}[theorem]{Lemma}
\newtheorem{proposition}[theorem]{Proposition}
\newtheorem{claim}[theorem]{Claim}
\theoremstyle{definition}
\newtheorem{definition}[theorem]{Definition}
\theoremstyle{remark}
\newtheorem{remark}[theorem]{Remark}

\def\ar{\mathfrak{A}[R]}
\def\arr{\mathfrak{A}_r[R]}
\def\kck{K\rtimes K^\times}
\def\rx{R^\times}
\def\kx{K^\times}
\def\cq{C_{\mathbb{Q}}}
\def\qn{\mc{Q}_{\mathbb{N}}}
\def\qz{\mc{Q}_{\mathbb{Z}}}
\def\cpk{C_{\mathrm{p}}^*(\kck,\mc{R})}

\newcommand{\mc}[1]{\text{$\mathcal{#1}$}}
\newcommand{\newf}[3]{{#1}:{#2}\longrightarrow {#3}}
\newcommand{\newfd}[5]{\begin{eqnarray*}{#1}:{#2} &\longrightarrow &{#3} \\{#4} &\longmapsto &{#5}\end{eqnarray*}}

\begin{document}

\title[Partial Crossed Product Description of the $C^*$-algebras of Domains]{Partial Crossed Product Description of the $C^*$-algebras Associated with Integral Domains}

\author[Giuliano Boava]{Giuliano Boava${}^1$}
\address{Giuliano Boava, Instituto Nacional de Matem\'atica Pura e Aplicada, 22460-320 - Rio de Janeiro/RJ, Brazil.}
\email{gboava@gmail.com}
\thanks{${}^1$Research supported by CNPq, Brazil.}

\author[Ruy Exel]{Ruy Exel${}^2$}
\address{Ruy Exel, Departamento de Matem\'atica, Universidade Federal de Santa Catarina, 88040-900 - Florian\'opolis/SC, Brazil.}
\email{r@exel.com.br}
\thanks{${}^2$Research partially supported by CNPq, Brazil.}

\begin{abstract}
Recently, Cuntz and Li introduced the $C^*$-algebra $\ar$ associated to an integral domain $R$ with finite quotients. In this paper, we show that $\ar$ is a partial group algebra of the group $\kck$ with suitable relations, where $K$ is the field of fractions of $R$. We identify the spectrum of this relations and we show that it is homeomorphic to the profinite completion of $R$. By using partial crossed product theory, we reconstruct some results proved by Cuntz and Li. Among them, we prove that $\ar$ is simple by showing that the action is topologically free and minimal.
\end{abstract}
\date{05 October 2010}
\maketitle

\section{Introduction}\label{introduction}

Fifteen years ago, motivated by the work of Julia \cite{julia}, Bost and Connes constructed a $C^*$-dynamical system having the Riemann $\zeta$-function as partition function \cite{bost-connes}. The $C^*$-algebra of the Bost-Connes system, denoted by $\cq$, is a Hecke $C^*$-algebra obtained from the inclusion of the integers into the rational numbers. In \cite{laca-raeburn2}, Laca and Raeburn showed that $\cq$ can be realized as a semigroup crossed product and, in \cite{laca-raeburn3}, they characterized the primitive ideal space of $\cq$.

In \cite{arledge-laca-raeburn}, \cite{cohen} and \cite{laca}, by observing that the construction of $\cq$ is based on the inclusion of the integers into the rational numbers, Arledge, Cohen, Laca and Raeburn generalized the construction of Bost and Connes. They replaced the field $\mathbb{Q}$ by an algebraic number field $K$ and $\mathbb{Z}$ by the ring of integers of $K$. Many of the results obtained for $\cq$ were generalized to arbitrary algebraic number fields (at least when the ideal class group of the field is $h=1$) \cite{laca-frankenhuijsen}, \cite{laca-larsen-neshveyev}.

Recently, a new construction appeared. In \cite{cuntz}, Cuntz defined two new $C^*$-algebras: $\qn$ and $\qz$. Both algebras are simple and purely infinite and $\qn$ can be seen as a $C^*$-subalgebra of $\qz$. These algebras encode the additive and multiplicative structure of the semiring $\mathbb{N}$ and of the ring $\mathbb{Z}$. Cuntz showed that the algebra $\qn$ is, essentially, the algebra generated by $\cq$ and one unitary operator. In \cite{yamashita}, Yamashita realized $\qn$ as the $C^*$-algebra of a topological higher-rank graph.

The next step was given by Cuntz and Li. In \cite{cuntz-li}, they generalized the construction of $\qz$ by replacing $\mathbb{Z}$ by a unital commutative ring $R$ (which is an integral domain with finite quotients by principal ideals). This algebra was called $\ar$. Cuntz and Li showed that $\ar$ is simple and purely infinite (when $R$ is not a field) and they related a $C^*$-subalgebra of its with the generalized Bost-Connes systems (when $R$ is the ring of integers in an algebraic number field having $h=1$ and, at most, one real place). In \cite{li}, Li extended the construction of $\ar$ to an arbitrary unital ring.

The aim of this text is to show that the algebra $\ar$ can be seen as a partial crossed product (when $R$ is an integral domain with finite quotients). We show that $\ar$ is isomorphic to a partial group algebra of the group $\kck$ with suitable relations, where $K$ is the field of fractions of $R$. By using the relationship between partial group algebras and partial crossed products, we see that $\ar$ is a partial crossed product by the group $\kck$. We characterize the spectrum of the commutative algebra arising in the crossed product and show that this spectrum is homeomorphic to $\hat{R}$ (the profinite completion of $R$). Furthermore, we show that the partial action is topologically free and minimal. By using that the group $\kck$ is amenable, we conclude that $\ar$ is simple.

Recently, some similar results appeared. In \cite{laca-raeburn4} and \cite{brownlowe-huef-laca-raeburn}, Brownlowe, an Huef, Laca and Raeburn showed that $\qn$ is a partial crossed product by using a boundary quotient of the Toeplitz (or Wiener-Hopf) algebra of the quasi-lattice ordered group $(\mathbb{Q}\rtimes \mathbb{Q}^\times_+, \mathbb{N}\rtimes \mathbb{N}^\times)$ (see \cite{nica} and \cite{laca-raeburn} for Toeplitz algebras of quasi-lattice ordered groups). We observe that our techniques are different from theirs. We don't use Nica's construction \cite{nica} (indeed, our group $\kck$ is not a quasi-lattice, in general). From our results, in the case $R=\mathbb{Z}$, we see that $\qz$ is a partial crossed product by the group $\mathbb{Q}\rtimes \mathbb{Q}^\times$. From this, it's imediate that $\qn$ is a partial crossed product by $\mathbb{Q}\rtimes \mathbb{Q}^\times_+$ (as in \cite{brownlowe-huef-laca-raeburn}).

Before we go to the main result we give, in the section \ref{preliminaries}, a breifly review about the algebra $\ar$ and the theories of partial crossed products and partial group algebras. In the section \ref{pgadescriptionAR}, we state our main theorem: the algebra $\ar$ is isomorphic to a partial group algebra. In the section \ref{pcpdescriptionAR}, we study $\ar$ by using the techniques of partial crossed products. We recover the faithful conditional expectation constructed by Cuntz and Li in \cite[Proposition 1]{cuntz-li} in a very natural way. Futhermore, we use the concepts of topological freeness and minimality of a partial action to show that $\ar$ is simple.

\section{Preliminaries}\label{preliminaries}

\subsection{The $C^*$-algebra $\ar$ of an Integral Domain}\label{adomain}

Throughout this text, $R$ will be an integral domain (unital commutative ring without zero divisors) with the property that the quotient $R/(m)$ is finite, for all $m\neq0$ in $R$. We denote by $\rx$ the set $R\backslash\{0\}$ and by $R^*$ the set of units in $R$.

\begin{definition}\cite[Definition 1]{cuntz-li}\label{definitionAR}
The {\bf regular $C^*$-algebra of} $R$, denoted by $\ar$, is the universal C*-algebra generated by isometries $\{s_m \ | \ m\in\rx\}$ and unitaries $\{u^n \ | \ n\in R\}$ subject to the relations
\begin{enumerate}[(CL1)]
\item $s_ms_{m'}=s_{mm'}$;
\item $u^nu^{n'}=u^{n+n'}$;
\item $s_mu^n=u^{mn}s_m$;
\item $\displaystyle\sum_{l+(m)\in R/(m)}u^ls_m s_m^*u^{-l}=1$;
\end{enumerate}
for all $m,m'\in\rx$ and $n,n'\in R$.
\end{definition}

We denote by $e_m$ the range projection of $s_m$, namely $e_m=s_ms_m^*$. It's easily seen that, under (CL2) and (CL3), $u^le_mu^{-l}=u^{l'}e_mu^{-l'}$ if $l+(m)=l'+(m)$. From this, we see that the sum in (CL4) is independent of the choice of $l$.

Let $\{\xi_r \ | \ r\in R\}$ be the canonical basis of the Hilbert space $\ell^2(R)$ and consider the operators $S_m$ and $U^n$ on $\ell^2(R)$ given by $S_m(\xi_r)=\xi_{mr}$ and $U^n(\xi_r)=\xi^{n+r}$.

\begin{definition}\cite[Section 2]{cuntz-li}\label{definitionARR}
The {\bf reduced regular $C^*$-algebra of} $R$, denoted by $\arr$, is the $C^*$-subalgebra of $\mc{B}(\ell^2(R))$ generated by the operators $\{S_m \ | \ m\in\rx\}$ and $\{U^n \ | \ n\in R\}$.
\end{definition}

One can checks that $S_m$ is an isometry, $U^n$ is a unitary and satisfy (CL1)-(CL4). Hence, there exists a surjective $*$-homomorphism $\ar\longrightarrow\arr$.

In \cite{cuntz-li}, Cuntz and Li showed that, when $R$ is not a field, $\ar$ is simple; therefore the above $*$-homomorphism is a $*$-isomorphism. In the section \ref{pcpdescriptionAR}, we will show that $\ar$ is simple (when $R$ is not a field) by using the partial crossed product description of $\ar$.

For future references, we need the following lemma, proved by Cuntz and Li:

\begin{lemma}\cite[Lemma 1]{cuntz-li}\label{lemmapcommute}
For all $n,n'\in R$ and $m,m'\in\rx$, the projections (in $\ar$) $u^ne_mu^{-n}$ and $u^{n'}e_{m'}u^{-n'}$ commute.
\end{lemma}

More details about these algebras can be found in \cite{cuntz}, \cite{cuntz-li}, \cite{cuntz-li2}, \cite{cuntz-li3}, \cite{larsen-li}, \cite{li} and \cite{yamashita}.

\subsection{Partial Crossed Products}\label{pcprodutcts}

Here, we review some basic facts about partial actions and partial crossed products.

\begin{definition}\cite[Definition 1.1]{dokuchaev-exel}\label{definitionPA}
A {\bf partial action} $\alpha$ of a (discrete) group $G$ on a $C^*$-algebra $\mc{A}$ is a collection $(\mc{D}_g)_{g\in G}$ of ideals of $\mc{A}$ and $*$-isomorphisms $\newf{\alpha_g}{\mc{D}_{g^{-1}}}{\mc{D}_g}$ such that
\begin{enumerate}[(P{A}1)]
\item $\mc{D}_e=\mc{A}$, where $e$ represents the identity element of $G$;
\item $\alpha_h^{-1}(\mc{D}_h\cap \mc{D}_{g^{-1}})\subseteq\mc{D}_{(gh)^{-1}}$;
\item $\alpha_g\circ \alpha_h(x)=\alpha_{gh}(x), \ \ \forall \: x\in\alpha_h^{-1}(\mc{D}_h\cap \mc{D}_{g^{-1}})$.
\end{enumerate}
\end{definition}

In the above definition, if we replace the $C^*$-algebra $\mc{A}$ by a locally compact space $X$, the ideals $\mc{D}_g$ by open sets $X_g$ and the $*$-isomorphisms $\alpha_g$ by homeomorphisms $\newf{\theta_g}{X_{g^{-1}}}{X_g}$, we obtain a {\bf partial action} $\theta$ of the group $G$ on the space $X$. A partial action $\theta$ on a space $X$ induces naturally a partial action $\alpha$ on the $C^*$-algebra $C_0(X)$. The ideals $D_g$ are $C_0(X_g)$ and $\alpha_g(f)=f\circ\theta_{g^{-1}}$.

We say that a partial action $\theta$ on a space $X$ is {\bf topologically free} if, for all $g\in G\backslash\{e\}$, the set $F_g=\{x\in X_{g^{-1}} \ | \ \theta_g(x)=x\}$ has empty interior. A subset $V$ of $X$ is {\bf invariant} under the partial action $\theta$ if $\theta_g(V\cap X_{g^{-1}})\subseteq V$, for every $g\in G$. The partial action $\theta$ is {\bf minimal} if there are no invariant open subsets of $X$ other than $\emptyset$ and $X$. It's easy to see that $\theta$ is minimal if, and only if, every $x\in X$ has dense orbit, namely $\mc{O}_x=\{\theta_g(x) \ | \ g\in G \ \mathrm{for \ which} \ x\in X_{g^{-1}}\}$ is dense in $X$.

\begin{definition}\cite[Definition 6.1]{dokuchaev-exel}\label{definitionPR}
A {\bf partial representation} $\pi$ of a (discrete) group $G$ into a unital $C^*$-algebra $\mc{B}$ is a map $\newf{\pi}{G}{\mc{B}}$ such that, for all $g,h\in G$,
\begin{enumerate}[(PR1)]
\item $\pi(e)=1$;
\item $\pi(g^{-1})=\pi(g)^*$;
\item $\pi(g)\pi(h)\pi(h^{-1})=\pi(gh)\pi(h^{-1})$.
\end{enumerate}
\end{definition}

From a partial action $\alpha$, we can construct two {\bf partial crossed products}: $\mc{A}\rtimes_{\alpha}G$ (full) and $\mc{A}\rtimes_{\alpha,r}G$ (reduced). We can define both as follows: let $\mc{L}$ be the normed $*$-algebra of the finite formal sums $\sum_{g\in G} a_g\delta_g$, where $a_g\in\mc{D}_g$. The operations and the norm in $\mc{L}$ are given by $(a_g\delta_g)(a_h\delta_h)=\alpha_g(\alpha_{g^{-1}}(a_g)a_h)\delta_{gh}$, $(a_g\delta_g)^*=\alpha_{g^{-1}}(a_g^*)\delta_{g^{-1}}$ and $||\sum_{g\in G} a_g\delta_g||=\sum_{g\in G} ||a_g||$. If we denote by $B_g$ the vector subspace $\mc{D}_g\delta_g$ of $\mc{L}$, then the family $(B_g)_{g\in G}$ generates a Fell bundle. The full and the reduced crossed products are, respectively, the full and the reduced cross sectional algebra of $(B_g)_{g\in G}$. It's well known that $\mc{A}\rtimes_{\alpha}G$ is universal with respect to a covariant pair $(\varphi,\pi)$, where $\newf{\varphi}{\mc{A}}{\mc{B}}$ is a $*$-homomorphism ($\mc{B}$ is a unital $C^*$-algebra), $\newf{\pi}{G}{\mc{B}}$ is a partial representation of $G$ and the covariant equations are $\varphi(\alpha_g(x))=\pi(g)\varphi(x)\pi(g^{-1})$ for $x\in\mc{D}_{g^{-1}}$ and $\varphi(x)\pi(g)\pi(g^{-1})=\pi(g)\pi(g^{-1})\varphi(x)$ for $x\in\mc{A}$.

There exists a faithful conditional expectation $\newf{E}{\mc{A}\rtimes_{\alpha,r}G}{\mc{A}}$ given by $E(a\delta_g)=a$ if $g=e$, and $E(a\delta_g)=0$ if $g\neq e$. When the Fell bundle $(B_g)_{g\in G}$ is amenable ($G$ amenable implies its), the full and reduced constructions are isomorphic and, in this case, there exists a faithful conditional expectation of $\mc{A}\rtimes_{\alpha}G$ onto $\mc{A}$.

There is a close relation between topological freeness and minimality of the partial action and ideals of the reduced crossed product. If $\theta$ is a topologically free partial action on a space $X$ then $\theta$ is minimal if, and only if, $C_0(X)\rtimes_{\alpha,r}G$ is simple, where $\alpha$ is the action induced by $\theta$. Under the amenability hypothesis, this result is valid for the full crossed product too.

For more details about partial crossed products, see \cite{dokuchaev-exel}, \cite{exel}, \cite{exel2}, \cite{exel3} and \cite{exel-laca-quigg}.

\subsection{Partial Group Algebras}\label{pgalgebras}

Let $G$ be a discrete group, let $\mc{G}$ be the set $G$ without the group operations and denote the elements in $\mc{G}$ by $[g]$ (namely, $\mc{G}=\{[g] \ | \ g\in G\}$). The {\bf partial group algebra of} $G$, denoted by $C_{\mathrm{p}}^*(G)$, is defined to be the universal $C^*$-algebra generated by the set $\mc{G}$ with the relations
$$\mc{R}_{\mathrm{p}}=\{[e]=1\}\cup \{[g^{-1}]=[g]^*\}_{g\in G}\cup \{[g][h][h^{-1}]=[gh][h^{-1}]\}_{g,h\in G}.$$
The algebra $C_{\mathrm{p}}^*(G)$ is universal with respect to a partial representation. Observe that the relations in $\mc{R}_{\mathrm{p}}$ correspond to the partial representation axioms (PR1), (PR2) and (PR3). Sometimes, we will refer to a relation in $\mc{R}_{\mathrm{p}}$ by indicating the corresponding axiom.

Consider the natural bijection between $\mc{P}(G)$ and $\{0,1\}^G$, where $\mc{P}(G)$ is the power set of $G$. With the product topology, $\{0,1\}^G$ is a compact Hausdorff space. Give to $\mc{P}(G)$ the topology of $\{0,1\}^G$. Denote by $X_G$ the subset of $\mc{P}(G)$ of the subsets $\xi$ of $G$ such that $e\in \xi$. Clearly, with the induced topology of $\mc{P}(G)$, $X_G$ is a compact space. For each $g\in G$, let $X_g=\{\xi\in X_G \ | \ g\in\xi\}$. It's easy to see that $\newf{\theta_g}{X_{g^{-1}}}{X_g}$ given by $\theta_g(\xi)=g\xi$ is a homeomorphism. The collection of open sets $(X_g)_{g\in G}$ of $X_G$ with the homeomorphisms $\theta_g$ define a partial action $\theta$ of $G$ on $X_G$. The partial crossed product $C(X_G)\rtimes_{\alpha}G$ is isomorphic to $C_{\mathrm{p}}^*(G)$ (where $\alpha$ is the partial action induced by $\theta$).

For each $g\in G$, we abreviate $[g][g^{-1}]$ by $e_g$. Let $\mc{R}$ be a set of relations on $\mc{G}$ such that every relation is of the form
$$\sum_i\lambda_i\prod_je_{g_{ij}}=0.$$
The {\bf partial group algebra of $G$ with relations} $\mc{R}$, denoted by $C_{\mathrm{p}}^*(G,\mc{R})$, is defined to be the universal $C^*$-algebra generated by the set $\mc{G}$ with the relations $\mc{R}_{\mathrm{p}}\cup\mc{R}$. Given a partial representation $\pi$ of $G$, we can extend $\pi$ naturally to sums of produtcs of elements in $\mc{G}$. If this extension satisfies the relations $\mc{R}$, we say that $\pi$ is a {\bf partial representation that satisfies} $\mc{R}$. The algebra $C_{\mathrm{p}}^*(G,\mc{R})$ is universal with respect to a partial representation that satisfies the relations $\mc{R}$.

Denote by $1_g$ the function in $C(X_G)$ given by $1_g(\xi)=1$ if $g\in\xi$ and $1_g(\xi)=0$ otherwise. By an abuse of notation, we also denote by $\mc{R}$ the subset of $C(X_G)$ given by the functions $\sum_i\lambda_i\prod_j1_{g_{ij}}$, where $\sum_i\lambda_i\prod_je_{g_{ij}}=0$ is a relation in (the original) $\mc{R}$. The {\bf spectrum of the relations} $\mc{R}$ is defined to be the compact Hausdorff space
$$\Omega_{\mc{R}} = \{\xi\in X_G\,|\,f(g^{-1}\xi)=0, \ \forall\,f\in\mc{R}, \ \forall\,g\in \xi\}.$$
Let $\Omega_g = \{\xi\in \Omega_{\mc{R}} \ | \ g\in\xi\}$. By restricting the above $\theta_g$ to $\Omega_{g^{-1}}$, we obtain a partial action (again denoted by $\theta$) of $G$ on $\Omega_{\mc{R}}$ (the open sets are the $\Omega_g$'s and the homeomorphisms are the restrictions of the $\theta_g$'s). The main result concerning $C_{\mathrm{p}}^*(G,\mc{R})$ says that this algebra is isomorphic to the partial crossed product $C(\Omega_{\mc{R}})\rtimes_{\alpha}G$ (again, $\alpha$ is the partial action induced by $\theta$).

The above results are proved in \cite{exel3} and \cite{exel-laca-quigg}.

\section{Partial Group Algebra Description of $\ar$}\label{pgadescriptionAR}

Let $R$ be an integral domain satisfying the conditions stated in the previous section. Denote by $K$ the field of fractions of $R$ and consider the semidirect product $\kck$. The elements of $\kck$ will be denoted by a pair $(u,w)$, where $u\in K$ and $w\in K^{\times}$. Recall that $(u,w)(u',w')=(u+u'w,ww')$ and $(u,w)^{-1}=(-u/w,1/w)$. We denote by $[u,w]$ an element of set $\kck$ without the group operations (as the set $\mc{G}$ associated to $G$ in the previous section).\footnote{Sometimes, we work with the element $(u,w)^{-1}$ or the element $(u_1,w_1)(u_2,w_2)$. For these elements, our corresponding notations will be $[(u,w)^{-1}]$ and $[(u_1,w_1)(u_2,w_2)]$.} Again, denote $[g][g^{-1}]$ by $e_g$. Consider the sets of relations
$$\mc{R}_1=\left\{e_{(n,1)}=1 \ | \ n\in R\right\}, \ \ \ \mc{R}_2=\left\{e_{\left(0,\frac{1}{m}\right)}=1 \ \bigl| \ m\in \rx\right\},$$
$$\mc{R}_3=\left\{\sum_{n+(m)\in R/(m)}e_{(n,m)}=1 \ \biggl| \ m\in\rx\right\}$$
and $\mc{R}=\mc{R}_1\cup\mc{R}_2\cup\mc{R}_3$. We observe that, under the relations $\mc{R}_1$ and $\mc{R}_{\mathrm{p}}$ (relations stated in the previous section), the sum in $\mc{R}_3$ independs of the choice of $n$. Indeed, for $k\in R$,
$$e_{(n+km,m)} = [n+km,m][(n+km,m)^{-1}] \stackrel{\mc{R}_1}{=} [(n,m)(k,1)]e_{(-k,1)}[(k,1)^{-1}(n,m)^{-1}] =$$
$$[(n,m)(k,1)][(k,1)^{-1}][k,1][(k,1)^{-1}(n,m)^{-1}] \stackrel{\mathrm{(PR3)}}{=}$$
$$ [n,m][k,1][(k,1)^{-1}][k,1][(k,1)^{-1}][(n,m)^{-1}] = [n,m]e_{(k,1)}e_{(k,1)}[(n,m)^{-1}] = e_{(n,m)}.$$

\begin{remark}\label{remarkR1}
The relations in $\mc{R}_1$ are unnecessary. They can be obtained from $\mc{R}_3$ with $m=1$.
\end{remark}

Consider the partial group algebra $\cpk$. We will show that this algebra is isomorphic to $\ar$.

\begin{proposition}\label{propositionPsi}
There exists a $*$-homomorphism $\newf{\Psi}{\ar}{\cpk}$ such that $\Psi(u^n)=[n,1]$ and $\Psi(s_m)=[0,m]$.
\end{proposition}

\begin{proof}
We need to show that $[n,1]$ is a unitary (for $n\in R$), that $[0,m]$ is an isometry (for $m\in\rx$) and that the relations (CL1)-(CL4) are satisfied. From $\mc{R}_1$ and (PR2), we have $[n,1][n,1]^*=e_{(n,1)}=1$ and $[n,1]^*[n,1]=e_{(-n,1)}=1$, ie, $[n,1]$ is a unitary. Similarly, from $\mc{R}_2$ and (PR2) we see that $[0,m]$ is an isometry. By using this fact,
$$\Psi(s_ms_{m'})=[0,m][0,m']=[0,m][0,m'][0,m']^*[0,m']\stackrel{\mathrm{(PR3)}}{=}$$
$$[0,mm'][0,m']^*[0,m']=[0,mm']=\Psi(s_{mm'}),$$
hence (CL1) is satisfied. We can prove (CL2) in the same way. To show (CL3), note that
$$\Psi(s_mu^n)=[0,m][n,1]=[0,m][n,1][n,1]^*[n,1]\stackrel{\mathrm{(PR3)}}{=}[mn,m][n,1]^*[n,1]=[mn,m],$$
because $[n,1]$ is a unitary. On the other hand,
$$\Psi(u^{mn}s_m)=[mn,1][0,m]=[mn,1][mn,1]^*[mn,1][0,m]\stackrel{\mathrm{(PR3)}}{=}$$
$$[mn,1][mn,1]^*[mn,m]=[mn,m].$$
Finally, (CL4) follows from $\mc{R}_3$ and\footnote{Be careful with the $e$'s! The notation $e_m$ represents $s_ms_m^*$ in $\ar$ and $e_{(n,m)}$ represents $[n,m][n,m]^*$ in $\cpk$.}
$$\Psi(u^ne_mu^{-n})=[n,1][0,m][0,m]^*[-n,1]=[n,m][0,1/m][-n,1][-n,1]^*[-n,1]\stackrel{\mathrm{(PR3)}}{=}$$
$$[n,m][(n,m)^{-1}][-n,1]^*[-n,1]=[n,m][(n,m)^{-1}]=e_{(n,m)}.$$
\end{proof}

Now, we will construct an inverse for $\Psi$. In the next claim, note that every element in $\kck$ can be writen under the form $\left(\frac{n}{m'},\frac{m}{m'}\right)$, where $n\in R$ and $m,m'\in\rx$.

\begin{claim}\label{claimPi}
The map $\newf{\pi}{\kck}{\ar}$ given by $\pi\left(\left(\frac{n}{m'},\frac{m}{m'}\right)\right)=s_{m'}^*u^ns_m$ is independent of the representation of $\left(\frac{n}{m'},\frac{m}{m'}\right)$.
\end{claim}

\begin{proof}
Let $\left(\frac{n}{m'},\frac{m}{m'}\right)=\left(\frac{q}{p'},\frac{p}{p'}\right)$, ie, $pm'=p'm$ and $m'q=p'n$. Hence,
$$s_{p'}^*u^qs_p = s_{p'}^*s_{m'}^*s_{m'}u^qs_p \stackrel{\mathrm{(CL3)}}{=} s_{p'}^*s_{m'}^*u^{m'q}s_{m'}s_p \stackrel{\mathrm{(CL1)}}{=} s_{p'm'}^*u^{m'q}s_{m'p} \stackrel{\mathrm{(CL1)}}{=}$$
$$s_{m'}^*s_{p'}^*u^{np'}s_{p'}s_m \stackrel{\mathrm{(CL3)}}{=} s_{m'}^*s_{p'}^*s_{p'}u^ns_m = s_{m'}^*u^ns_m.$$
\end{proof}

\begin{proposition}\label{propositionPi}
The map $\pi$ defined above is a partial representation of $\kck$ that satisfies $\mc{R}$.
\end{proposition}

\begin{proof}
First, we will show that $\pi$ is a partial representation. Since $\pi((0,1))=s_1^*u^0s_1=1$, we have (PR1). Observe that
$$\pi\left(\left(\frac{n}{m'},\frac{m}{m'}\right)^{-1}\right) = \pi\left(\left(\frac{-n}{m},\frac{m'}{m}\right)\right) = s_m^*u^{-n}s_{m'} = \pi\left(\left(\frac{n}{m'},\frac{m}{m'}\right)\right)^*,$$
which shows (PR2). To see (PR3), let $g=\left(\frac{q}{p'},\frac{p}{p'}\right)$ and $h=\left(\frac{n}{m'},\frac{m}{m'}\right)$. We have $gh=\left(\frac{m'q+pn}{p'm'}, \frac{pm}{p'm'}\right)$ and, therefore,
$$\pi(gh)\pi(h^{-1}) = \pi(gh)\pi(h)^* = (s_{p'm'}^*u^{m'q+pn}s_{pm})(s_m^*u^{-n}s_{m'}) \stackrel{\mathrm{(CL1), (CL2), (CL3)}}{=}$$
$$s_{p'}^*u^qs_{m'}^*s_pu^ns_ms_m^*u^{-n}s_{m'} = s_{p'}^*u^qs_{m'}^*s_p\underbrace{u^ns_ms_m^*u^{-n}}_{}\underbrace{s_{m'}s_{m'}^*}_{}s_{m'} \stackrel{\mathrm{Lemma} \ \ref{lemmapcommute}}{=}$$
$$s_{p'}^*u^qs_{m'}^*s_ps_{m'}s_{m'}^*u^ns_ms_m^*u^{-n}s_{m'} \stackrel{\mathrm{(CL1)}}{=} (s_{p'}^*u^qs_p)(s_{m'}^*u^ns_m)(s_m^*u^{-n}s_{m'}) = \pi(g)\pi(h)\pi(h^{-1}).$$
This shows that $\pi$ is a partial representation. It remains to show that the extention of $\pi$ satisfies the relations in $\mc{R}$. By remark \ref{remarkR1}, it suffices to show that the relations in $\mc{R}_2$ and $\mc{R}_{3}$ are satisfied. It follows from
$$\pi(e_{(0,1/m)})=\pi([0,1/m][0,m])=s_m^*u^0s_1s_1^*u^{0}s_m=1$$
and
$$\pi\left(\sum_{n+(m)\in R/(m)}e_{(n,m)}\right)=\sum_{n+(m)\in R/(m)}s_1^*u^ns_ms_m^*u^{-n}s_1 = 1.$$
\end{proof}

\begin{remark}\label{remarkgeneralPi}
We can define $\pi$ for a general representation of a element in $\kck$ by $\pi\left(\left(\frac{n}{m''},\frac{m}{m'}\right)\right)=s_{m''}^*u^ns_{m'}^*s_{m''}s_m$.
\end{remark}

\begin{theorem}\label{theoremisomorphism}
The $*$-homomorphism $\Psi$ defined above is a $*$-isomorphism. Its inverse $\newf{\Phi}{\cpk}{\ar}$ is given by $\Phi\left(\left[\frac{n}{m'},\frac{m}{m'}\right]\right) = s_{m'}^*u^ns_m$.
\end{theorem}

\begin{proof}
The existence of $\Phi$ follows from $\pi$ and the universal property of $\cpk$. It remains to show that $\Psi$ and $\Phi$ are inverses each other. Indeed, $\Phi(\Psi(u^n))=\Phi([n,1])=s_1^*u^ns_1=u^n$, $\Phi(\Psi(s_m))=\Phi([0,m])=s_1^*u^0s_m=s_m$ and
$$\Psi\left(\Phi\left(\left[\frac{n}{m'},\frac{m}{m'}\right]\right)\right) = \Psi(s_{m'}^*u^ns_m) = \left[0,1/m'\right]\left[n,1\right]\left[0,m\right] = $$
$$ \left[0,1/m'\right]\left[0,1/m'\right]^*\left[0,1/m'\right]\left[n,1\right]\left[n,1\right]^*\left[n,1\right]\left[0,m\right] = \left[\frac{n}{m'},\frac{m}{m'}\right].$$
\end{proof}

\section{Partial Crossed Product Description of $\ar$}\label{pcpdescriptionAR}

Before characterizing $\ar$ as a partial crossed product, note that the group $\kck$ is solvable and, hence, amenable. Therefore, there exists a faithful conditional expectation (imported from the partial crossed product realization) $\newf{E}{\cpk}{C^*(\{e_g\}_{g\in\kck})}$ given by
$$E([g_1][g_2]\cdots[g_k])=\delta_{g_1g_2\cdots g_k,e}[g_1][g_2]\cdots[g_k].$$
In \cite[Proposition 1]{cuntz-li}, Cuntz and Li constructed a faithful conditional expectation $\Theta$ on $\ar$ given by $\Theta(s_{m''}^*u^ns_ms_m^*u^{-n'}s_{m'})=\delta_{m',m''}\delta_{n,n'}s_{m'}^*u^ns_ms_m^*u^{-n}s_{m'}$. The next proposition shows that, under the $*$-isomorphism $\Psi$, $E$ and $\Theta$ are the same conditional expectation.

\begin{proposition}\label{propositionCE}
$E\circ\Psi=\Psi\circ\Theta$.
\end{proposition}

\begin{proof}
First of all, observe that $\left(\frac{n}{m''},\frac{m}{m''}\right)\left(\frac{-n'}{m},\frac{m'}{m}\right)=(0,1)$ if, and only if, $m'=m''$ and $n=n'$. Hence,
$$E\circ\Psi(s_{m''}^*u^ns_ms_m^*u^{-n'}s_{m'}) = E\left(\left[\frac{n}{m''},\frac{m}{m''}\right]\left[\frac{-n'}{m},\frac{m'}{m}\right]\right) =$$
$$\delta_{m',m''}\delta_{n,n'}\left[\frac{n}{m'},\frac{m}{m'}\right]\left[\frac{-n}{m},\frac{m'}{m}\right].$$
On the other hand
$$\Psi\circ\Theta(s_{m''}^*u^ns_ms_m^*u^{-n'}s_{m'}) = \Psi(\delta_{m',m''}\delta_{n,n'}s_{m'}^*u^ns_ms_m^*u^{-n}s_{m'}) =$$
$$\delta_{m',m''}\delta_{n,n'}\left[\frac{n}{m'},\frac{m}{m'}\right]\left[\frac{-n}{m},\frac{m'}{m}\right].$$
\end{proof}

We already know that $\ar$ is a partial crossed product. Indeed, every partial group algebra is a partial crossed product (see section \ref{pgalgebras}). From now on, our goal is to study $\ar$ by this way.

There exists a natural partial order on $\rx$ given by the divisibility: we say that $m\leq m'$ if there exists $r\in R$ such that $m'=mr$. Whenever $m\leq m'$, we can consider the canonical projection $\newf{p_{m,m'}}{R/(m')}{R/(m)}$. Since $(\rx, \leq)$ is a directed set, we can consider the inverse limit
$$\hat{R} = \lim_{\longleftarrow}\{R/(m), \ p_{m,m'}\},$$
which is the {\bf profinite completion of} $R$. In this text, we shall use the following concrete description of $\hat{R}$:
$$\hat{R}=\left\{(r_m+(m))_m\in\prod_{m\in\rx}R/(m) \ \biggl| \ p_{m,m'}(r_{m'}+(m'))=r_m+(m), \ \mathrm{if} \ m\leq m' \right\}.$$
Give to $R/(m)$ the discrete topology, to $\prod_{m\in\rx}R/(m)$ the product topology and to $\hat{R}$ the induced topology of $\prod_{m\in\rx}R/(m)$. With the operations defined componentwise, $\hat{R}$ is a compact topological ring. There exists a canonical inclusion of $R$ into $\hat{R}$ given by $r\longmapsto(r+(m))_m$ (to see injectivity, take $r\neq0$, $m$ non-invertible and note that $r\notin(rm)$).

The above partial order can be extended to $\kx$. For $w,w'\in\kx$, we say that $w\leq w'$ if there exists $r\in R$ such that $w'=wr$. Denote by $(w)$ the fractional ideal generated by $w$, namely $(w)=wR\subseteq K$. As before, if $w\leq w'$, we can consider the canonical projection\footnote{By the second isomorphism theorem, it could be $\newf{p_{w,w'}}{R/(R\cap(w'))}{R/(R\cap(w))}$.} $\newf{p_{w,w'}}{(R+(w'))/(w')}{(R+(w))/(w)}$. As before, we consider the inverse limit
$$\hat{R}_K = \lim_{\longleftarrow}\{(R+(w))/(w), \ p_{w,w'}\}\cong$$
$$\left\{(u_w+(w))_w\in\prod_{w\in\kx}(R+(w))/(w) \ \biggl| \ p_{w,w'}(u_{w'}+(w'))=u_w+(w), \ \mathrm{if} \ w\leq w' \right\}.$$
It is a compact topological ring too. In fact, $\hat{R}_K$ is naturally isomorphic to $\hat{R}$ as topological ring. In this text, we use $\hat{R}_K$ instead of $\hat{R}$ to simplify our proofs.

It's easy to see that, when $R$ is a field, then $\hat{R}\cong\hat{R}_K\cong\{0\}$.

Let $\Omega$ be the spectrum of the relations $\mc{R}$ (see section \ref{pgalgebras}). We will show that $\Omega$ is homeomorphic to $\hat{R}_K$ (hence, homeomorphic to $\hat{R}$). Define
\newfd{\rho}{\hat{R}_K}{\mc{P}(\kck)}{(u_w+(w))_w}{\{(u_w+rw,w) \ | \ w\in\kx, \ r\in R\}.}
Note that the definition is independent of the choice of $u_w$ in $u_w+(w)$.

\begin{claim}\label{claimRho}
$\rho(\hat{R}_K)\subseteq\Omega$.
\end{claim}

\begin{proof}
Let $(u_w+(w))_w\in\hat{R}_K$. By the definition of $\hat{R}_K$, if $w\leq w'$, then $u_{w'}=u_{w}+kw$ for some $k\in R$. Denote $\rho((u_w+(w)))$ by $\xi$. Clearly, $(0,1)\in\xi$. We need to show that $f(g^{-1}\xi)=0$, for all $f\in\mc{R}$ and $g\in\xi$. Fix $g=(u_w+rw,w)\in\xi$. Let $f=1_{(n,1)}-1$ in $\mc{R}_1$ and note that $f(g^{-1}\xi)=0$ is equivalent to $g(n,1)\in\xi$. Since $g(n,1)=(u_w+rw,w)(n,1)=(u_w+(r+n)w,w)$, we have $g(n,1)\in\xi$. Now, let $f=1_{(0,1/m)}-1$ in $\mc{R}_2$. Similarly, we must show that $g(0,1/m)\in\xi$. Observe that $g(0,1/m) = (u_w+rw,w)(0,1/m) = (u_w+rw,w/m)$. Since $w/m\leq w$, then $g(0,1/m) = (u_{w/m} + k(w/m) + rw,w/m)=(u_{w/m} + (k+rm)(w/m))\in\xi$. To finish, fix $m\in\rx$ and let $f=\sum_{n+(m)}1_{(n,m)}-1$ in $\mc{R}_3$. We must show that there exists one, and only one class $n+(m)$ such that $g(n,m)\in\xi$. Indeed, $g(n,m) = (u_w+rw,w)(n,m) = (u_w+(n+r)w,wm) = (u_{wm} + (n+r-k)w,wm)$ and, for it belongs to $\xi$, we must have $(n+r-k)w\in(wm)$. Hence, $n\equiv k-r\mod m$, in other words, there exists only one class $n+(m)$ such that $g(n,m)\in\xi$. Since $\mc{R}=\mc{R}_1\cup\mc{R}_2\cup\mc{R}_3$, the proof is completed.
\end{proof}

\begin{proposition}\label{propositionRho}
$\newf{\rho}{\hat{R}_K}{\Omega}$ is a homeomorphism.
\end{proposition}

\begin{proof}
\ \\
{\bf Injectivity.} Let $(u_w+(w))_w,(v_w+(w))_w\in\hat{R}_K$ such that $\rho((u_w+(w)))=\rho(v_w+(w)))$. By the definition of $\rho$, the elements in $\rho((u_w+(w)))$ whose second component equals $w$ are of the form $(u_w+rw,w)$. Since $(v_w,w)\in\rho((v_w+(w)))$ and, therefore, $(v_w,w)\in\rho((u_w+(w)))$, we must have $v_w=u_w+rw$ for some $r\in R$. This show that $(u_w+(w))_w=(v_w+(w))_w$. \\
{\bf Surjectivity.} Let $\xi\in\Omega$. The relations in $\mc{R}_1$ and $\mc{R}_2$ together implies that if $g\in\xi$, then $g(q/p,1/p)\in\xi$ for all $q\in R$ and $p\in\rx$ (fix $g$ and apply $f(g^{-1}\xi)=0$ for various $f$). For each $m\in\rx$, let $f=\sum_{n+(m)}1_{(n,m)}-1$ in $\mc{R}_3$ and apply $f(g^{-1}\xi)=0$ with $g=(0,1)$ to see that there exists only one class $n+(m)$ such that $(n,m)\in\xi$. Denote this class by $u_m+(m)$. Since $g(0,1/p)\in\xi$ if $g\in\xi$, then $p_{m,mp}(u_{mp}+(mp))=(u_{m}+(m))$. From this, we can define unambiguously $u_w+(w)=u_m+(w)$ for $w=m/m'\in\kx$. One can see that the classes $u_w+(w)$ are compatible with the projections $p_{w,w'}$ by using that $g(q/p,1/p)\in\xi$ if $g\in\xi$. Hence, we have constructed $(u_w+(w))_w\in\hat{R}_K$. We claim that $\rho((u_w+(w)))=\xi$. Since $(u_w,w)\in\xi$, $(u_w,w)(q,1)=(u_w+qw,w)$ must belongs to $\xi$. This shows that $\rho((u_w+(w)))\subseteq\xi$. Suppose, by contradiction, $\rho((u_w+(w)))\neq\xi$. Hence, there exists $h\in\xi$ such that $h\notin\rho((u_w+(w)))$. If we write $h=(n'/m',m/m')$, then $h\notin\rho((u_w+(w)))$ is equivalent to $n'-m'u_{m}\notin(m)$. Let $g=(u_m,1/m')$, $h'=(u_m,m/m')$ and note that both belong to $\rho((u_w+(w)))$ (hence, belong to $\xi$). Since $g^{-1}h=(-m'u_m,m')(n'/m',m/m')=(n'-m'u_m,m)$, $g^{-1}h'=(0,m)$ and $n'-m'u_{m}\notin(m)$, then $f(g^{-1}\xi)\neq0$ if $f=\sum_{n+(m)}1_{(n,m)}-1$, which contradicts the fact that $\xi\in\Omega$. Hence, $\rho((u_w+(w)))=\xi$. \\
To finish the proof, observe that $\hat{R}_K$ and $\Omega$ are compact Hausdorff, therefore it suffices to show that $\rho$ (or $\rho^{-1}$) is continuous to conclude that $\rho$ is a homeomorphism. We will prove that $\rho^{-1}$ is continuous by showing that $\pi_w\circ\rho^{-1}$ is continuous for all $w\in\kx$, where $\newf{\pi_w}{\hat{R}_K}{(R+(w))/(w)}$ is the canonical projection. Since $(R+(w))/(w)$ is discrete, it suffices to show that $\rho\circ\pi_w^{-1}(\{u_w+(w)\})$ is an open set of $\Omega$, for all $u_w+(w)\in (R+(w))/(w)$. To see this, note that
$$\rho\circ\pi_w^{-1}(\{u_w+(w)\})=\{\xi\in\Omega \ | \ (u_w,w)\in\xi\},$$
which is an open set of $\Omega$ (recall that the topology on $\Omega$ is induced by the product topology of $\{0,1\}^{\kck}$).
\end{proof}

Following the section \ref{pgalgebras}, there exists a partial action of $\kck$ on $\Omega$. By the above proposition, we can define this partial action on $\hat{R}_K$. Let $\hat{R}_g=\rho^{-1}(\Omega_g)$, where $\Omega_g=\{\xi\in\Omega \ | \ g\in\xi\}$, and $\theta_g$ be the homeomorphism between $\hat{R}_{g^{-1}}$ and $\hat{R}_g$. It's easy to see that
$$\hat{R}_{(u,w)}=\{(u_{w'}+(w'))_{w'}\in\hat{R}_K \ | \ u_w + (w) = u + (w)\}$$
and
$$\theta_{(u,w)}((u_{w'}+(w'))_{w'}) = (u + wu_{w'} + (ww'))_{ww'} = (u + wu_{w^{-1}w'} + (w'))_{w'},$$
ie, $\theta_{(u,w)}$ acts on $\hat{R}_{(u,w)^{-1}}$ by the affine transformation corresponding to $(u,w)$. The next proposition, whose proof is trivial, will be useful later.

\begin{proposition}\label{propositionRg}
We have that
\begin{enumerate}[(i)]
\item $\hat{R}_{(u,w)} = \emptyset \ \ \Longleftrightarrow \ \ u\notin R+(w)$;
\item $\hat{R}_{(u,w)} = \hat{R}_K \ \ \Longleftrightarrow \ \ R\subseteq u+(w)$.
\end{enumerate}
\end{proposition}

Now, we describe the topology on $\hat{R}_K$. Since $\hat{R}_K$ is a singleton set when $R$ is a field, we shall assume that $R$ is not a field in this paragraph. For $w\in\kx$ and $C_w\subseteq (R+(w))/(w)$, we define the open set
$$V_w^{C_w}=\{(u_{w'}+(w'))_{w'}\in\hat{R}_K \ | \ u_w + (w)\in C_w\}.$$
Clearly, if $w\leq w'$, then $V_w^{C_w}=V_{w'}^{C_{w'}}$, where $C_{w'}=\{u+(w')\in(R+(w'))/(w') \ | \ u+(w)\in C_w\}$. From the product topology, we know that the finite intersections of open sets $V_w^{C_w}$ form a basis for the topology on $\hat{R}_K$. By taking a commom multiple of the $w$'s in the intersection, we see that every basic open set is of the form $V_w^{C_w}$ (since $V_w^{C_1}\cap V_w^{C_2}=V_w^{C_1\cap C_2}$). Futhermore, if $C_w\neq\emptyset$, $r$ is a non-invertible element in $R$ and $V_w^{C_w}=V_{wr}^{C_{wr}}$, then $C_{wr}$ has, at least, two elements. Indeed, let $u+(w)\in C_w$ and $r_1,r_2\in R$ such that $r_1+(r)\neq r_2+(r)$. It's easy to see that $u+wr_1+(wr)$ and $u+wr_2+(wr)$ are in $C_{wr}$ and that $u+wr_1+(wr)\neq u+wr_2+(wr)$. This says that, if $V_w^{C_w}$ is non-empty, we can suppose that $C_w$ has more than one element.

\begin{proposition}\label{propositiontopfree}
The partial action $\theta$ on $\hat{R}_K$ is topologically free if, and only if, $R$ is not a field.
\end{proposition}

\begin{proof}
If $R$ is a field, then $\hat{R}_K=\{0\}$ and, hence, $\theta$ is not topologically free. Conversely, suppose that $R$ is not a field. We need to show that $F_g=\{x\in \hat{R}_{g^{-1}} \ | \ \theta_g(x)=x\}$ has empty interior, for all $g\in\kck\backslash\{(0,1)\}$. We shall consider two cases: $g=(u,1)$ and $g=(u,w)$, $w\neq1$. \\
{\bf Case 1.} If $u\notin R$, then the proposition \ref{propositionRg} says that $\hat{R}_{g^{-1}}=\emptyset$. So, we can suppose $u\in R$. If $F_g\neq\emptyset$, then equation $\theta_g(x)=x$ implies that $u\in (m)$ for every $m\in\rx$. Since $R$ is not a field, then $u=0$. This show that $F_g=\emptyset$ if $g=(u,1)$ and $u\neq0$. \\
{\bf Case 2.} Let $g=(u,w)$ such that $w\neq1$ and $u\in R+(w)$ (if $u\notin R+(w)$, then $\hat{R}_{g^{-1}}=\emptyset$). Let $V$ be a non-empty open set contained in $\hat{R}_{g^{-1}}$. We will show that there exists $x\in V$ such that $\theta_g(x)\neq x$. By shrinking $V$ if necessary, we can suppose that $V=V_{w'}^{C_{w'}}$. Futhermore, we can assume that $C_{w'}$ has more than one element. Let $u_1+(w')$ and $u_2+(w')$ be distinct elements of $C_{w'}$, hence $u_1-u_2\notin(w')$. Suppose, by contradiction, $\theta_g(x)= x$ for all $x\in V$. Since $(u_i+(w''))_{w''}\in V$, $i=1,2$, then
$$\theta_{(u,w)}((u_i+(w''))_{w''})=(u_i+(w''))_{w''} \ \ \Longrightarrow \ \ (u+wu_i + (w''))_{w''}=(u_i+(w''))_{w''}.$$
By choosing $w''=(w-1)w'$ (note that $w\neq1$), we see that $u+(w-1)u_i\in((w-1)w')$, for $i=1,2$. By subtracting the equations (for different $i$'s), we have $(w-1)(u_1-u_2)\in((w-1)w')$ and, therefore $u_1-u_2\in(w')$; which is a contradiction! This show that $F_g$ has empty interior.
\end{proof}

\begin{proposition}\label{propositionminimal}
The partial action $\theta$ is minimal.
\end{proposition}

\begin{proof}
If $R$ is a field, then the result is trivial. Now, suppose that $R$ is not a field. We will prove that every $x\in\hat{R}_K$ has dense orbit (see section \ref{pcprodutcts}) by showing that if $V$ is a non-empty open set, then there exists $g\in\kck$ such that $x\in\hat{R}_{g^{-1}}$ and $\theta_g(x)\in V$. Let $x=(u_w+(w))_w\in\hat{R}_K$ and $V=V_{w'}^{C_{w'}}$ non-empty. Take $u'+(w')\in C_{w'}$ and observe that we can suppose, without loss of generality, $u'\in R$ and $u_{w'}\in R$. Let $g=(u'-u_{w'},1)$. By the proposition \ref{propositionRg}, $\hat{R}_{g^{-1}}=\hat{R}_K$ and, hence, $x\in\hat{R}_{g^{-1}}$. To finish, note that $\theta_g(x)=\theta_{(u'-u_{w'},1)}((u_w+(w))_w) = (u'-u_{w'}+u_{w}+(w))_w\in V$.
\end{proof}

Following, we summarize the results of this section.

\begin{theorem}\label{theoremPCProduct}
The algebra $\ar$ is $*$-isomorphic to the partial crossed product $C(\hat{R}_K)\rtimes_\alpha\kck$, where $\alpha$ is the partial action induced by $\theta$. The $*$-isomorphism is given by $u^n\longmapsto1\delta_{(n,1)}$ and $s_m\longmapsto1_{(0,m)}\delta_{(0,m)}$, where $1_{(0,m)}$ is the characteristic function of $\hat{R}_g$.
\end{theorem}

\begin{theorem}\label{theoremsimple}
$\ar$ is simple.
\end{theorem}

\begin{proof}
By the propositions \ref{propositiontopfree} and \ref{propositionminimal}, the reduced crossed product $C(\hat{R}_K)\rtimes_{\alpha,r}\kck$ is simple. Since $\kck$ is amenable, then $C(\hat{R}_K)\rtimes_\alpha\kck\cong C(\hat{R}_K)\rtimes_{\alpha,r}\kck$ and, therefore, $C(\hat{R}_K)\rtimes_\alpha\kck$ is simple. The result follows from the previous theorem.
\end{proof}

\begin{corollary}\label{corollaryArArr}
$\ar\cong\arr$.
\end{corollary}

When $R=\mathbb{Z}$, we can restrict our partial action to the subgroup $\mathbb{Q}\rtimes\mathbb{Q}_+^*$ of $\mathbb{Q}\rtimes\mathbb{Q}^*$ and the corresponding partial crossed product is the algebra $\qn$ introduced by Cuntz in \cite{cuntz} and realized as a partial crossed product in \cite{brownlowe-huef-laca-raeburn} by Brownlowe, an Huef, Laca and Raeburn.

\end{document}